\documentclass[a4paper,11pt]{article}

\usepackage{amsmath}
\usepackage{amsthm}
\usepackage{amssymb}
\usepackage{cancel}

\usepackage{epsfig}
\usepackage{psfrag}
\usepackage{subfigure}
\usepackage{graphicx}

\usepackage[mathscr,mathcal]{eucal}
\usepackage{enumerate}

\usepackage{t1enc}
\usepackage[latin2]{inputenc}

\def \N {\mathbb N}

\def \Z {\mathbb Z}

\def\cB{\mathcal{B}}

\def\vareps{\varepsilon}
\def \eps {\epsilon}

\def\chn {{\breve n}}

\newcommand{\prob}[1]{\ensuremath{\mathbf{P}\left(#1\right)}}

\newcommand{\condprob}[2]{\ensuremath{\mathbf{P}\big(#1\bigm|#2\big)}}

\newcommand{\abs}[1]{\left|\,{#1}\,\right|}
\newcommand{\norm}[1]{\left\|\,{#1}\,\right\|}

\def \wt {\widetilde}
\def\ol{\overline}
\def\ul{\underline}
\def\wh{\widehat}

\newtheorem {theorem}{Theorem}

\newtheorem {lemma}[theorem]{Lemma}
\newtheorem {corollary}[theorem]{Corollary}
\newtheorem {proposition}[theorem]{Proposition}

\def\be{\begin{equation}}
\def\ee{\end{equation}}

\def\bea{\begin{eqnarray}}
\def\eea{\end{eqnarray}}

\title{Stuck Walks}

\author{
{\sc Anna Erschler}
\qquad
{\sc B\'alint T\'oth}
\qquad
{\sc Wendelin Werner}
}

\date {}

\begin{document}

\maketitle

\begin{abstract}
We investigate the asymptotic behaviour of a class of self-interacting nearest neighbour random walks on the one-dimensional integer lattice which are pushed by a particular linear combination of their own local time on edges in the neighbourhood of their current position. We prove that in a range of the relevant parameter of the model such random walkers can be eventually confined to a finite interval of length depending on the parameter value. The phenomenon arises as a result of competing self-attracting and self-repelling effects where in the named parameter range the former wins.



\end{abstract}

\section{Introduction and main result}
\label{s:intro}

Let $(X_n, n \ge 0)$ be a nearest neighbour path on the one-dimensional integer lattice $\Z$, and define for each $n\in\N$ and  $j\in\Z$, its local time $\ell (n,j)$
 on unoriented edges:
$$
\ell(n,j):=\#\{1\le m\le n\,:\, \{X_{m-1}, X_m \}=\{j-1,j\}\}.
$$
Throughout this paper the unoriented edge connecting the sites $j-1$ and $j$ will be denoted by $j$.

We fix a real parameter $\alpha$ and define
\begin{align}
\label{Delta}
\Delta(n,j)
:=
-\alpha\ell(n,j-1)+\ell(n,j)-\ell(n,j+1)+\alpha\ell(n,j+2)
\end{align}
for all $j \in \Z$ and $n \ge 0$. We then also define
$ \Delta_n  = \Delta (n, X_n)$, which is therefore
a particular linear combination of the number of visits by $X$ before time $n$ to the edges near $X_n$.

We consider a special type of self-interacting random walk $(X_n, n \ge 0)$ with long memory started from $X_0=0$,
whose law is described by the following ``dynamics'': for all $n \ge 0$,
\begin{align}
\label{law}
\condprob{X_{n+1} = X_n \pm1}{{\mathcal F}_n} =
\frac{\exp\{\pm\beta\Delta_n\}} {\exp\{\beta\Delta_n\}+\exp\{-\beta\Delta_n\}}
\end{align}
where $\beta>0$ is another fixed parameter of the problem and ${\mathcal F}_n = \sigma (X_0, \ldots, X_n)$.
In plain words, if $ \Delta_n$ is positive (respectively, negative), then the walker will prefer to jump to the right (resp., to the left) at its $(n+1)$-st jump.

 We are interested in the long time asymptotic behaviour of the walk. The parameter $\alpha$ plays
a  crucial role. Depending on its value the qualitative behaviour varies spectacularly. The role of the parameter $\beta$ is less dramatic.

The driving mechanism \eqref{law} is a generalization of the rules governing the so-called ``true'' self-avoiding random walk (or true self-repelling walk -- we will refer to is as the TSRW) in 1d.
Choosing $\alpha=0$ we obtain the TSRW with edge repulsion (the latter looks at each step at the number of times it has previously jumped along its two neighbouring edges, and favors the less-visited one)
 while choosing $\alpha=-1$ corresponds to the TSRW with site repulsion. In these two cases non-degenerate scaling limits for $n^{-2/3}X(n)$ are proved \cite{toth_95}, \cite{toth_werner_98}, respectively conjectured \cite{amit_parisi_peliti_83}, \cite{obukhov_peliti_83}, \cite{peliti_pietronero_87}. See the survey \cite{toth_01} for more information about these two cases.

It turns out that depending on the value of the parameter $\alpha$, a rather rich phase diagram emerges. We refer to \cite {erschler_toth_werner_10} for background and motivation.
For $\alpha$ in $[-1,1/3)$, we expect a similar scaling behaviour as for the TSRW. It is intuitively clear that for $\alpha$ close to $0$ the mechanism can be viewed as a minor perturbation of the TSRW case. The fact that the interval of parameters where this kind of asymptotic behaviour is expected is exactly $\alpha\in[-1,1/3)$ follows from more detailed arguments relying on the fact that this is the range of parameters where the coefficients of the linear combination defining $\Delta$ in \eqref{Delta} correspond to a positive definite sequence. For details of this argument see \cite{erschler_toth_werner_10}, \cite{ttv}. For $\alpha\in(-\infty,-1)$, when the walk is repelled by its past visits to its neighbouring edges and
even more strongly by its second-neighbouring edges,
a kind of slowing down phenomenon seems to occur, where the walk gets slowed down by self-built trapping environments. A more detailed discussion can be found in \cite {erschler_toth_werner_10}.

The results of the present paper will concern the range of values where $\alpha$ is positive.
In this case, the walk is repelled by its local time on the edges adjacent to its current position, but attracted by its previous visits to its two next-to-neighbouring edges. As we shall see, when $\alpha>1/3$, the self-attractiveness can win and the walk can remain stuck forever on a finite interval of consecutive sites, while this fails to hold for $\alpha\le1/3$. It is therefore natural to define
the (possibly infinite) interval ${\mathcal L}$ of points that are visited infinitely often by the walk.

We will use  a simple explicit sequence of values $(\alpha_{L}, L \ge 1)$ that decays to $1/3$ defined as follows: $\alpha_{1} = + \infty$ and for all $L \ge 2$,
$$ 
\alpha_{L} = \frac {1} { 1 + 2 \cos (2 \pi/(L+2))}.
$$
Our main result is the following:
\begin {theorem}
\label {thm:main}
 Suppose that $L \ge 1$. Then:
\begin  {itemize}
 \item
 If $\alpha \in (\alpha_{L+1}, \alpha_{L})$, then the probability that $\# {\mathcal L}= L +2$ is positive.
\item
 If $\alpha < \alpha_{L+1}$, then almost surely, $\# {\mathcal L} > L +2 $.
\end {itemize}
\end {theorem}

Note our convention of denoting discrete interval length: When $\# {\mathcal L} = L+2$, this means that the number of \emph{interior} lattice sites of ${\mathcal L}$ is $L$. Thus, $L+1$ will be the number of lattice edges in the interval (i.e., the length of the interval) and $L+2$ will be the number of sites in the closed interval, including the endpoints. Such a discrete interval will be of the type $\{x,x+1, \dots,x+L,x+L+1\}$, $x\in\Z$, with endpoints $x$ and $x+L+1$.

\begin {figure}[htbp]
\begin {center}
\includegraphics [height=1.9in]{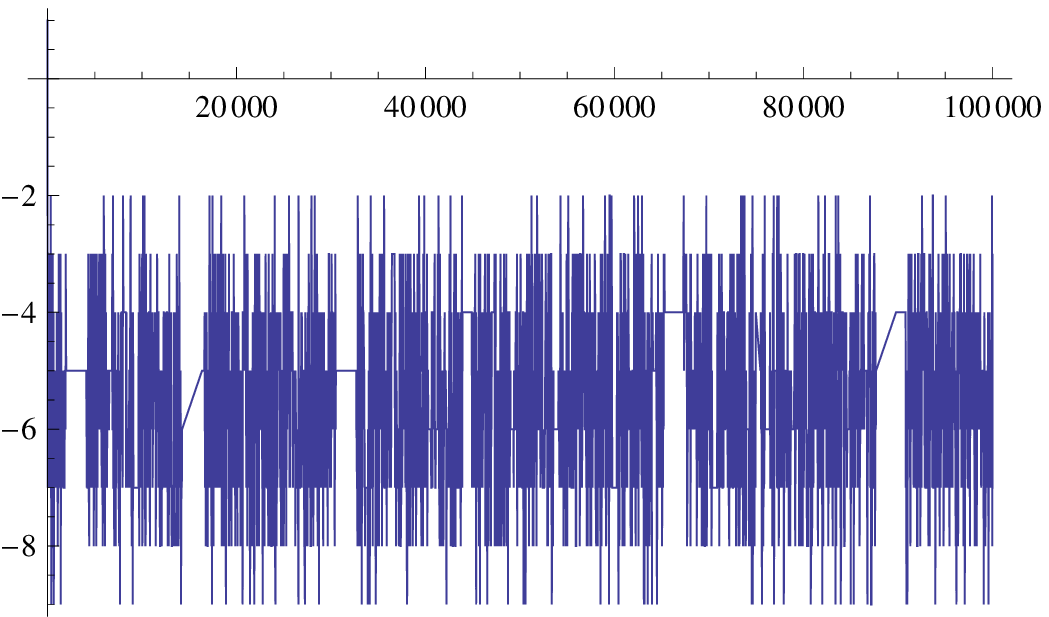}
\includegraphics [height=1.9in]{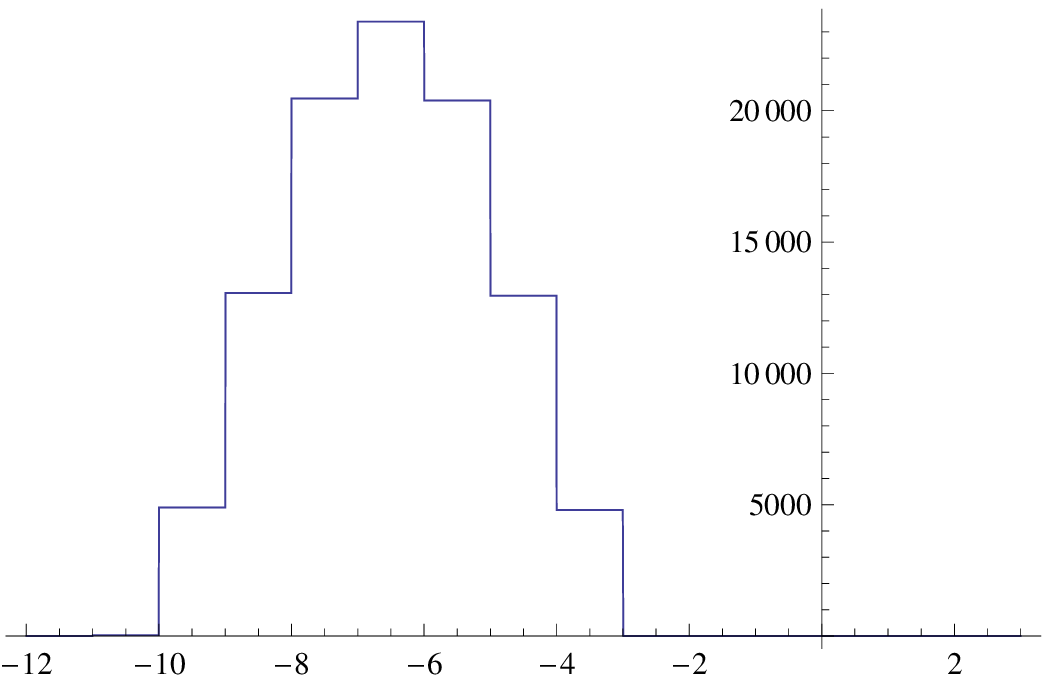}
\caption {A trajectory and its local time when $\alpha = 0.4$}
\end {center}
\end {figure}

It implies in particular that $\# {\mathcal L} \in \{ 0, \infty \}$ almost surely when $\alpha \le 1/3 = \inf_{L \ge 2} \alpha_{L} $ i.e. that the range of the walk is infinite.
This trapping phenomenon when $\alpha > 1/3$ is reminiscent of the asymptotic behaviour of the \emph{vertex reinforced random walk} in 1d, cf. \cite{pemantle_92}, \cite{pemantle_volkov_99}, \cite{tarres_04}. However, the differences are also conspicuous. The reinforcement scheme is of a  different type, and here, there is no clear monotone attractiveness in the self-interaction mechanism, rather a competition between self-attraction and self-repulsion, where the self-attraction wins. Due to this, the size of the trapping range increases to infinity as the parameter value approaches the borderline between the confined and non-confined regimes, at $\alpha =1/3$.

We will in fact also give a more precise description of the asymptotic behavior of the walk, or rather of its local times on edges, in the case where it is trapped. In the scenario that we will describe and happens with positive probability, the renormalized local time profile $ (n^{-1} \ell (n, j ), j \in \Z)$ will become deterministic
in the large-scale limit: We will describe for each $\alpha$ in $(\alpha_{L+1}, \alpha_L)$ an explicit sequence $u_1, \ldots, u_{L}$ (that will follow the values of a $\sin^2$ curve along an arithmetic sequence) so that when
${\mathcal L} = \{  x , \ldots , {x+L+1} \}$ for some $x\in\Z$, one then necessary has (up to a set of zero probability)
$$  \lim_{n \to \infty} \frac {1}{n} ( \ell (n,x+1), \ldots, \ell (n,x+ L+1)) =  (u_1, \ldots, u_{L+1}).$$

\section{Probabilistic part of the proof}
\label{s:forcibly_confined_walk}

Let $L \ge 1$ and $\alpha < \alpha_{L}$ be fixed throughout most of this section (note that by definition $\alpha_{1}  = \infty$).
We define an auxiliary nearest neighbour random walk $(Y_n , n \ge 0)$ confined \emph{by force} to the interval $[0,L+1]$.
We will study some properties of $Y$ and see whether it is possible to couple $X$ and $Y$ is such a way that (with positive probability) they coincide forever.

By abuse of notation we will also denote by $\ell(n,j)$, respectively, by  $\Delta(n,j)$ and ${\mathcal F}_n$ the local time on edges, respectively, the linear combinations and the $\sigma$-field defined for $Y$ just as for $X$. The law of this walk started from $Y_0= 0$ is described by its dynamics
\begin {align*}
\condprob{Y_{n+1}= Y_n \pm1}{{\mathcal F}_n}
=
\left\{
\begin{array}{cl}
\displaystyle
\frac{1\pm1}{2}
& \text{ if } Y_n=0,\\ [13pt]
\displaystyle
\frac{\exp\{\pm\beta\Delta_n\}} {\exp\{\beta\Delta_n\}+\exp\{-\beta\Delta_n\}}
&\text{ if } Y_n \in\{1,\dots,L\},\\ [13pt]
\displaystyle
\frac{1\mp1}{2}
& \text{ if } Y_n =L+1.
\end{array}
\right.
\end{align*}
So, $Y$ behaves exactly as $X$ except that when it is on the boundary of the interval $[0, L+1]$, it is forced to jump inwards.

When $Y_n \notin \{0, L+1 \}$, we can interpret $\Delta_n$ as a local stream felt by the walker (due to its past) at time $n$.
If $\Delta_n$ is positive, it will tend to jump to the right, whereas when $\Delta_n$ is negative, it will tend to jump to the left.
We say that when it jumps (from $Y_n$ to $Y_{n+1}$) to the \emph{opposite direction} than the one suggested by the sign of $\Delta_n$ and when $Y_n \in [1, L]$, it does an upstream jump of intensity $|\Delta_n|$. We will be interested  in the relation between the maximal stream that the walker has ``successfully overcome'' before time $n$ and the maximal value of $\Delta (n,j)$.

The main ingredient of the proof of Theorem \ref {thm:main} is the following \emph{deterministic} statement, that says that 
there is no way of building up a stream larger than $D \ge D_0$ somewhere in the interior of the interval without having earlier performed an upstream jump of intensity larger $\vareps D$ somewhere.
Its proof is given in section \ref{s:proof_of_propo_mds}.

\begin{proposition}
\label{prop:mds}
Suppose that $\alpha < \alpha_{L}$.
There exist constants $\vareps=\vareps(\alpha, L)>0$ and $D_0=D_0(\alpha, L)<\infty$ such that
for any nearest neighbour walk trajectory $(Y_n, n \ge 0) $ in $\{0,\dots,L+1\}$, any $D \ge D_0$ and any $n$, 
at least one of the following two statements hold:
\begin {itemize}
 \item 
For all $j\in\{1,\dots,L\}$, and a positive $n$ such that  $\abs{\Delta(n,j)}\le D$.
\item 
During its first $n$ steps, the walk $Y$ has performed at least one upstream jump of intensity larger than
$\vareps D$.
\end {itemize}
\end{proposition}

In other words, 
\medbreak

Note that for any $n \ge 0$, if $\Delta_n >0$ and $Y_n \not= L+1$, the conditional probability that $Y_{n+1} - Y_n = -1$ given ${\mathcal F}_n$ is smaller than
$\exp ( - 2 \beta \Delta_n)$. The symmetric result holds when $Y_n \not= 0$. It follows readily that, for any positive $n$ and $D$, the probability that $Y$ does an upstream jump of intensity greater than $\varepsilon D$ at its $n$-th jump is smaller than $\exp (-2 \beta \varepsilon D)$.
Hence, the proposition implies that for all $D > D_0$ and all $n \ge 0$,
$$ \prob {\max_{j\in\{1,\dots,L\}} \abs{\Delta(n,j)} \ge D }
\le n e^{-2 \beta \vareps D}.
$$
A Borel-Cantelli argument immediately implies that for the walk $Y$:

\begin{corollary}
\label{co:nobigstream}
There exists a constant $c$, such that almost surely,
$ \max_{j\in\{1,\dots,L\}} \abs{\Delta(n,j)} \le c \log n$
for all large $n$.
\end{corollary}

We see in particular that when $n$ is very large, all $L$ values
$\Delta (n,1), \ldots, \Delta(n,L)$ are very small compared to $n$. We can keep in mind that these are
simple linear combinations of the $L+1$ non-negative numbers $\ell (n,1), \ldots, \ell (n, L+1)$ that also satisfy
$$ \ell (n,1) + \ldots + \ell (n, L+1) = n . $$

This leads us naturally to study the set of possible solutions $(l_1, \ldots, l_{L+1})$ to the linear system of $L+1$ equations
given by
\begin {equation}
\label {refu}
d_1 = d_2 = \cdots = d_L = 0 \hbox { and } l_1 + \cdots + l_{L+1} = 1,
\end {equation}
where $d_j = - \alpha l_{j-1} + l_j - l_{j+1} + \alpha l_{j+2}$, with the convention $l_{0}=l_{L+2} = 0$.

Note that the conditions $d_{1} = \ldots = d_L = 0$ mean that $ l_0=0 , \ldots, l_{L+1}, l_{L+2}=0$ are part of a bi-infinite ``Fibonacci-type'' sequence
$(\tilde l_j,j \in \Z)$ that satisfies $\tilde d_j=0$ for all $j \in \Z$ (with obvious notation). Note that then,
$$d_0 =  - l_1 + \alpha l_2 = l_0 - l_1 + \alpha l_2 = \tilde d_0 + \alpha \tilde l_{-1} =  \alpha \tilde l_{-1}$$ and 
similarly $d_{L+1} = -\alpha \tilde l_{L+3}$.
Such bi-infinite sequences are on a periodic curve as soon as $\alpha > 1/3$ and then, if we define  $\omega \in (0, \pi)$ by
$$ \cos (\omega ) = \frac {1- \alpha}{2 \alpha}, $$
 such a bi-infinite sequence is necessarily of the form
$$ \tilde l_j = A + B \cos ( \omega  j + \varphi ) $$
for some $A$, $B$ and $\varphi$. The fact that $\tilde l_0 = 0 $ shows that it is possible to take
$$ \tilde l_j = B ( \cos (\omega j + \varphi ) - \cos (\varphi)).$$

It is then
immediate to check that when $\alpha < \alpha_{L}$,
one can take
$$
\varphi  (\alpha)
:=\frac {1}{2} ( {2\pi} - (L+1) \omega)
$$
and that  $(u_1, \ldots , u_{L+1})$ given by the formula
$$
u_j = \frac { \cos(\varphi ) - \cos(\omega j + \varphi )}{Z} = \frac { \sin^2 ( (\varphi /2) + j ( \omega/2) ) - \sin^2 ( \varphi/2)}{Z'}
$$
(where $Z$ and $Z'$ are the normalisation constants chosen so that $u_1 + \cdots + u_{L+1} = 1$) is the solution to our system of equations.

Similarly, when $\alpha \le 1/3$, the solution to (\ref {refu}) can be easily worked out. 
In any case, we can observe that it is non-negative (as soon as $\alpha < \alpha_L$), and that:
\begin {itemize}
 \item
When $\alpha  < \alpha_{L+1}$ and $(l) = (u)$,
$d_0<0$ and $d_{L+1} > 0$.
\item When
$\alpha \in ( \alpha_{L+1} , \alpha_{L} ) $ and $(l) = (u)$,
$d_0 >0$ and $d_{L+1} < 0$.
\end {itemize}
This will be important later.

Corollary \ref{co:nobigstream} therefore implies in particular that $\Delta (n, j) / n \to 0$ almost surely when $ n \to \infty$ for all $j \in \{1, \ldots,L \}$,
  which in turn implies that almost surely,
$$
\lim_{n \to \infty } \frac{\ell(n,j)}{n} = u_j
$$
for all $j=1, \ldots , L+1$ (this reasoning in fact implies an asymptotic upper bound on the rate of convergence).
Note that this implies also that almost surely,
$\Delta (n, 0) / n \to  d_0$ and $\Delta (n, L+1) / n \to  d_{L+1} $ when $n \to \infty$.

\medbreak

In order to prove Theorem \ref{thm:main}, we now have to see if it is possible to couple $X$ and $Y$ in such a way that they stick together with positive probability.
We first try to couple the walks starting at the time $0$. Given the similar dynamics, the optimal way to couple them is quite obvious: If we first define $(X_n)$, we can then simply define $Y_n =X_n$ as long as $ n \le \tau $, where
$$\tau  := \inf \{ t \ge 0 \ : \ X_t \notin [0, L+1] \}.$$
 In order for $\tau$ to be greater than $t$, it therefore suffices that at each of the times
$n \in \{ 0, 1, \ldots, t \}$ at which $Y_n \in \{ 0, L+1 \}$ (we call ${\mathcal N}_t$ this random set of times), $X$ jumps inwards (i.e., not out of our interval).
Hence,
$$
\condprob{ \tau >  t}{ Y_0, \ldots, Y_t }  = \prod_{n \in {\mathcal N}_t}    \left(\frac { e^{\beta \Delta_n} 1_{\{Y_n=0\}} + e^{-\beta \Delta_n}1_{\{Y_n=L+1\}} }{ e^{\Delta_n} + e^{-\Delta_n}}
\right).$$
The previous description of the asymptotic behavior of $Y$ (and of $\Delta ( n, L+1)$ and $\Delta (n, 0)$) immediately implies on the one hand that $\tau < \infty$ almost surely if $\alpha < \alpha_{L+1}$ (because $\Delta (n, 0)$ tends to $-\infty$), and that on the other hand, the probability that $\tau = \infty$ is strictly positive (by a simple Borel-Cantelli argument, due to the fact that $\Delta(n, 0) \sim d_0 n $ and $\Delta (n, L+1) \sim d_{L+1} n $ almost surely when $n \to \infty$) when $\alpha \in (\alpha_{L+1}, \alpha_{L})$.

\medbreak
Hence, we have proved that:
\begin {itemize}
\item When $\alpha < \alpha_{L+1}$, the probability that $X$ stays forever in $\{0, 1, \ldots, L+1\}$ is zero.
\item When $\alpha \in (\alpha_{L+1}, \alpha_{L})$, the probability that $X$ stays forever in $\{0, \ldots, L+1 \}$ is positive. Furthermore, in this case, the asymptotic local time profile of $Y$ (and therefore also of $X$) satisfies $\ell(j, n) / n \to u_j$ as $n \to \infty$ for all $j \in \{1, \ldots, L+1\}$.
\end {itemize}

To conclude the proof of Theorem \ref {thm:main}, it remains to notice that when $\alpha < \alpha_{L+1}$, the previous argument can be immediately adapted to show that for all $m\ge 1$ and for all finite nearest-neighbour sequence
$x_0, \ldots, x_{m}$, the conditional probability that for all $n \ge n_0$, $X_n \in \{x_{m}, 1+ x_{m}, \ldots , L+1+ x_m \}$ given $\{ X_0 = x_0, \ldots, X_{m}=x_m \}$ is zero. It suffices  to couple $(X_{m+n} - X_m, n \ge 0)$ with $Y$, where the local time has been suitably initialized.
Similarly, the ``successful'' coupling argument between $X$ and $Y$ can be started after these $m$ steps and shows that when $\alpha \in (\alpha_{L+1}, \alpha_{L})$ and the walk gets eventually trapped in $[x, x+ L+1]$ for some $x$, then
$(\ell (n,x+1), \ldots , \ell (n, x+L+1) )  \sim   (n u_1, \ldots, n u_{L+1})$ when $n \to \infty$.

\begin {figure}[htbp]
\begin {center}
\includegraphics [height=1.8in]{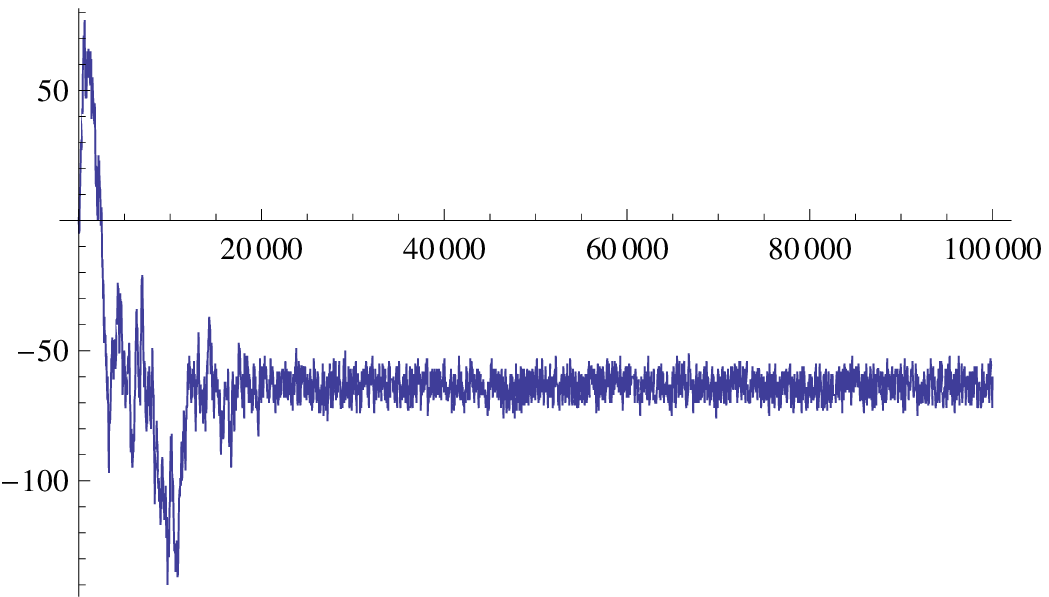}
\includegraphics [height=1.8in]{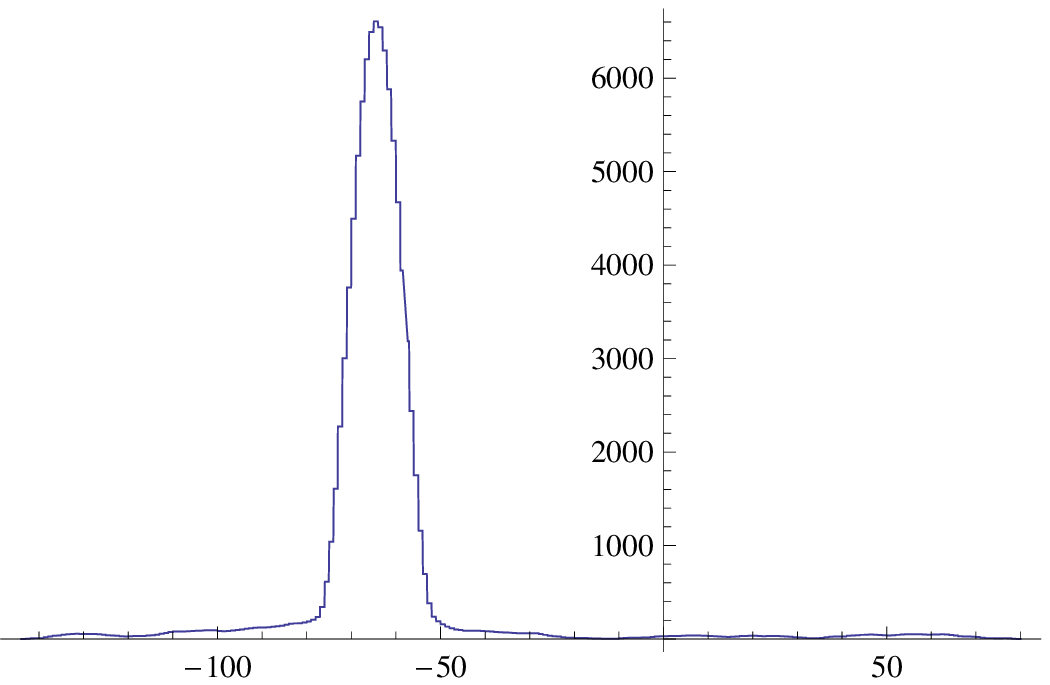}
\caption {The trapping phenomenon: A trajectory and its local time when $\alpha = 0.34$}
\end {center}
\end {figure}

\section{Proof of the main combinatorial statement}
\label{s:proof_of_propo_mds}

Our goal in this section is to prove Proposition \ref{prop:mds}. As we have already noticed, this is a deterministic statement about nearest-neighbour paths.

When $L=1$, the statement turns out to be straightforward (the walk is confined to two edges). So, we restrict ourselves to $L \ge 2$ (so that $\alpha < \alpha_{L} \le \alpha_2 =  1$).
Recall that $(u_1, \ldots, u_{L+1})$ is the unique solution to the linear system (\ref {refu}) (with unknown $(l_1, \ldots, l_{L+1})$)
and that the values of $u_j$ for $j =1, \ldots, L+1$ are all positive.
By continuity, it follows that by choosing $\gamma$ small enough, we can ensure that, if
\begin   {equation}
\label {system}
 \max (| d_1 |, \ldots, |d_L| ) \le \gamma \hbox { and } l_1 + \ldots + l_{L+1} = 1,
 \end {equation}
then the $l_j$'s are as close as we wish to the $u_j$'s. It follows (from the corresponding signs of $d_0$ and $d_{L+1}$ for $(u)$) that:
\begin {itemize}
 \item If $\alpha < \alpha_{L+1}$, there exists $\gamma ( \alpha, L)$ such that if  $\gamma \le \gamma (\alpha, L)$, then (\ref {system}) implies that
$d_0< 0$ and $d_{L+1}>0 $.
\item If $\alpha \in (\alpha_{L+1}, \alpha_{L})$, there exists $\gamma (\alpha, L)$ such that if   $\gamma \le  \gamma (\alpha, L)$, then (\ref {system}) implies that
$d_0 > 0$ and $d_{L+1} < 0$.
\end {itemize}

Let us now choose our constants $D_0$ and $\vareps$.
Recall that $\alpha < \alpha_{L} \le \alpha_{2} =1$. From now on, we choose (and fix) $\gamma$ such that
$$ \gamma < \frac {1}{100} \min (1, \gamma (\alpha, L), \gamma (\alpha, L-1), \ldots, \gamma (\alpha, 2)).$$
We then define
$$
\vareps:=\gamma^{L(L+1)} \hbox { and } D_0:= 10/ \vareps.
$$
The main role of $D_0$ will be to ensure that all time-intervals that we will talk about are non-empty (we will omit to mention this throughout the proofs) as soon as $D \ge D_0$.

Suppose that $n\mapsto Y_n \in\{0,\dots,L+1\}$ is a given nearest-neighbour trajectory.
First we define some particular times, corresponding to first appearances of streams of certain intensity.
For all  $M\in(0,\infty)$, we let (using the convention $\inf \emptyset = \infty$)
\begin{align*}
\theta_{+}(M)
&:=
\inf\{n:\Delta(n,j)\ge \gamma^{jL} M \text{ for some } j\in\{1,\dots,L\}\},
\\
\theta_{-}(M)
&:=
\inf\{n:\Delta(n,L+1-j)\le - \gamma^{jL} M \text{ for some } j\in\{1,\dots,L\}\},
\end {align*}
and we finally let $\sigma (M)$ denote the first time at which the walk makes an upstream jump of intensity
greater than $M$.

\medbreak
We will prove that with our choice of  $\vareps$ and $D_0$, then necessarily, $\theta_+ (D) \ge \sigma (\vareps D )$ for all $D \ge D_0$. By symmetry, we then
necessarily also have $\theta_- (D) \ge \sigma (\vareps D)$. But because
 $\gamma<1$, this implies that
$$
\sigma ( \vareps D) \le \min ( \theta_+ (D) , \theta_- (D) ) \le \inf\{n:\max_{j\in\{1,\dots,L\}} | \Delta(n,j)| \ge D \},
$$
which completes the proof of the proposition.

\medbreak

We will prove by contradiction that $\theta_+ (D) \ge \sigma (\vareps D )$. We will therefore from now on assume that $Y$ belongs to the set $\cB_+$ of paths $Y$ such that $\theta_+(D) < \sigma(\vareps D)$ for some given $D \ge D_0$ i.e. that $Y$ has created a ``very strong'' stream to the right (in the sense defined by $\theta_+(D)$) without having done any upstream jump of intensity larger than $\vareps D$. 

\medbreak

We start with two simple preliminary remarks.

\begin{lemma}
\label{lemma:conf}
(i)
For any $n_1<n_2$ and $j\in\{1,\dots,L\}$,
$$\abs{\Delta(n_2,j)- \Delta(n_1,j)}\le n_2-n_1.$$ 
(ii)
Let $1<M<\infty$ be fixed, $j \in\{1,\dots,L\}$ and $n_1<n_2<\sigma(M)$. If
\begin{align*}
\min_{n\in[n_1,n_2]}\Delta(n,j)> M+1
\end{align*}
then the walk is confined to the interval $\{j,\dots,L+1\}$ for the whole time-span $[n_1,n_2]$.
\end{lemma}

(ii) means that if there is a strong positive stream somewhere, then the walk is located to the right of this stream, unless it has performed a strong upstream jump before.
The symmetric result also holds, i.e., if there is a strong negative stream somewhere, then (unless the walk has performed a strong upstream jump before) the walk is to the
left of the stream.

\begin{proof}
(i)
This is straightforward, since for all $j\in\{1,\dots,L\}$ and $n\in\N$,
\begin{align}
\Delta(n+1,j)-\Delta(n,j)\in\{-1,-\alpha,0,\alpha, 1\}.
\end{align}
(ii)
Let
\begin{align*}
\hat n:=\max\{n<n_1:\Delta(n,j)\le M+1\}
.\end{align*}
Then $\Delta (\hat n + 1, j) > M+1 $, so that $\Delta (\hat n , j) > M$. Since $\Delta ( \hat n , j) < \Delta (\hat n +1 , j)$, one necessarily has
 $\{Y(\hat n), Y(\hat n +1)\}=\{j-1,j\}$ or $\{Y(\hat n), Y(\hat n +1)\}=\{j+1,j+2\}$.
But, the possibility $\{ Y(\hat n)=j,  Y(\hat n +1)=j-1 \}$ is excluded, since this would correspond to an upstream jump which is in conflict with the assumption $n_1<\sigma(M)$.
Hence, $Y( \hat n + 1) \ge j$ and  then, between $\hat n$ and $n_2$, the strong stream at $j$ will prevent the walk from jumping to the left of $j$.
\end{proof}

Since $\theta_+ (D)$ is finite, we can define the following:
\begin{align*}
\ol n & := \theta_+ (D)  = \min \{ n  \ : \ \Delta (n, j) \ge \gamma^{jL} D \hbox { for some } j \in \{1, \ldots, L \} \}
\\
J
&:=
\max\{j\in\{1,\dots,L\}: \Delta(\ol n,j)\ge \gamma^{jL} D\},
\\
\ul n
&:=
 \max\{n\le\ol n: \Delta(n,J)\le{\gamma^{JL} D}/ {2}\}.
\end{align*}
Throughout the proof we restrict ourselves to the time-span $[\ul n, \ol n]$ (this is the time-interval where we will detect a contradiction i.e. the necessity of an upstream jump -- in fact, we will zoom into  smaller time-intervals). 
Note that by the definition of $\ol n$ and $\ul n$ and (i) of Lemma \ref{lemma:conf},
\begin {equation}
\label {ref26}
\ol n - \ul n \ge {\gamma^{JL} D}/ {2}.
\end {equation}
Furthermore, the definition of $\ul n$ and our choice of $\gamma$, $\vareps$ and $D$ show that
$$
\min_{n\in[\ul n, \ol n]}\Delta(n,J)
= \Delta ( \ul n, J) \ge
\frac {\gamma^{JL} D}{2} - 1 
>
2\vareps D > \vareps D + 1,
$$
and thus, because of Lemma \ref {lemma:conf}-(ii),  any walk $Y$ in $\cB_+$ is confined to the interval $[J,L+1]$ for the whole time-span $[\ul n,\ol n]$ (recall that we assume that $\ol n \le \sigma (\eps D)$). It is important to notice here that the interval $[J, L+1]$ is strictly shorter than $[0,L+1]$.

Let us choose $\wh n$ to be the smallest integer such that
$$
\wh n  \ge \ul n + \frac{\gamma^{JL} D}{2}
$$
(mind that the symbol $\wh n$ will be used as a temporary variable which will be redefined in various ``subroutines'' of the proof) and note that, because of (\ref {ref26}), 
$ \wh n \le \ol n$. 
We now define
\begin{align*}
l_i
&:=
\ell(\wh n, J+i)-\ell(\ul n, J+i),
&&
i\in\{1,\dots,L-J+1\},
\\
d_i
&:=
\Delta(\wh n, J+i)-\Delta(\ul n, J+i),
&&
i\in\{1,\dots,L-J\}.
\end{align*}
Then $(l_1, \ldots, l_{L-J+1})$ is a family of non-negative numbers, and
$$ \| l \| := | l_1+ \ldots + l_{L-J+1} |  = \hat n - \ul n \ge
 \frac{\gamma^{JL} D}{2}.
$$
Let us define $$C = C(J,D) = \gamma^{L(J+1)} D.$$
 The reader might want to keep in mind that $\gamma$ is small and that 
$ D \ge C \ge \vareps D$.
We will soon prove the following lemma:
\begin{lemma}
\label{lemma:bound_on_Delta}
For all $Y \in \cB_+$, for all $j\in\{J+1,\dots,L\}$,
$$
\max_{n\in[\ul n,\ol n]}\abs{\Delta(n,j)}
\le
C \left(\frac{1}{\gamma}\right)^{j-(J+1)} .
$$
\end{lemma}
Let us now show how it implies the proposition. Note that then, for all $j \in \{ J+1, \ldots, L \}$,
$$
\max_{n\in[\ul n,\ol n]}  \abs{\Delta(n,j)}
\le
D \gamma^{(J+1)L} \gamma^{J+1 - L } \le D \gamma^{JL + 1 } \le 2 \gamma \| l \| .
$$ 
Hence
$$
\max_{j\in\{1,\dots,L -J +1 \}} \abs{d_j} \le 4 \gamma \norm{l}.
$$
Recall that $L+1-J \le L$; our choice of $\gamma \le \gamma (\alpha, L-J+1)$ therefore implies that
$$
d_0 = \Delta(\wh n, J)-\Delta(\ul n,J) < 0,
$$
which {is in contradiction} with the definition of $\ul n$. We conclude that $\cB_+$ is indeed empty, and that the proposition holds.
\medbreak

It now remains to prove Lemma \ref{lemma:bound_on_Delta}:
First we note that by the definition of $\ol n$, for all $j\in\{J+1,\dots,L\}$,
$$
\max_{n \le \ol n} \Delta(n,j)
\le
{\gamma^{jL}} D
\le \gamma^{(J+1)L} D = C \le
C \left(\frac{1}{\gamma}\right)^{j-(J+1)}.
$$
So, it remains to control the negative streams, i.e. to prove that for all $j \in \{ J+1, \ldots, L \}$,
\begin{align}
\label{remains}
\min_{n\in[\ul n,\ol n]} \Delta(n,j)
\ge
-C \left(\frac{1}{\gamma}\right)^{j-(J+1)} .
\end{align}
We will proceed by induction for $j=J+1,J+2, \dots,L$, from the left to the right.

\medbreak
We start with the case where $j=J+1$. 
We proceed in two steps. First we prove the following slightly stronger bound at time $\ul n$:
\begin{align}
\label{bound_beginning_1}
\Delta(\ul n, J+1)\ge-\frac {C}{2}.
\end{align}
Assume that the contrary holds i.e., that $ \Delta(\ul n, J+1) < - C/2$. Note that the definitions of $\gamma$ and $\vareps$ show that 
$$ \frac {C}{4}  =  \frac {\gamma^{L(J+1)} D}{4} \ge  \frac {\gamma^{L \times L} D}{4} \ge 2 \gamma^{L(L+1)} D =   2 \vareps D.$$
Furthermore, note that $C/4 \le \ol n - \ul n$.
Hence, from Lemma \ref{lemma:conf} and the fact that $\ol n \le \sigma (\vareps D)$,  it follows that $Y$ is confined in the interval  $[J,J+1]$ during $C/4$ steps,  
i.e, that it bounces back and forth on
this single edge for at least $C/{4}$ steps after $\ul n$. If we let
$\chn $ denote the integer value of $ \ul n + C/4$, it therefore follows that on the one hand $\chn \le \ol n$ and that on the other hand
$$
\Delta(\chn, J) - \Delta(\ul n, J)<0,
$$
which contradicts the definition of $\ul n$. Hence, we conclude that  \eqref{bound_beginning_1} indeed holds.

\medbreak
Next, we want to study what happens on the entire interval $[ \ul n, \ol n ]$. The quantity $\Delta (n, J+1)$ has to decrease from above $-C/2$ to below $-C$.
We let
\begin{align*}
&
\wh n:=\inf\{n\ge\ul n: \Delta(n,J+1)<-C\}
\\
&
\wt n:= \max\{n\le\wh n: \Delta(n,J+1)\ge -{C}/{2}\}
\end{align*}
(note that these definitions of  $\wh n$  and $\wt n$ are also temporary and will be redefined in other ``subroutines'' of the proof).
Assume that
$
\wh n \le \ol n$.
From \eqref{bound_beginning_1} it follows that (for all $Y \in \cB_+$),
$ \wt n>\ul n$ so that $[\wt n, \wh n ]\subset [\ul n , \ol n]$. But the definition of $\wt n$ shows that
\begin{align*}
\max_{n\in[\wt n, \wh n]} \Delta(n,J+1) \le -\frac{C}{2} + 1 < - \vareps D,
\end{align*}
so that we can deduce as before (using the fact that $\sigma (\vareps D) \ge \ol n$) from
 Lemma \ref{lemma:conf} that on the whole time span $[\wt n, \wh n]$ the walk is confined to the interval $[J, J+1]$ and thus we readily get
$$
\Delta(\wh n,J+1)-\Delta(\wt n, J+1) > 0,
$$
which contradicts the definition of $\wt n$ and $\wh n$. We conclude that (for all $Y \in \cB_+$) the bound \eqref{remains} holds for $j=J+1$.

\medbreak
The induction step follows next. If $J=L-1$ then we are done. So, assume that $J<L-1$, we let  $K\in\{J+1,\dots,L-1\}$ and assume that \eqref{remains} holds for $j\in\{J+1,\dots,K\}$. Our goal is to prove it for $j=K+1$. Again, first we divide this into two steps. Let us first prove the following slightly stronger bound at time $\ul n$:
\begin{align}
\label{bound_beginning_2}
\Delta(\ul n, K+1)
\ge
-\frac{C}{2} \left(\frac{1}{\gamma}\right)^{K-J} .
\end{align}
Assume the contrary. At time $\ul n$, we have right stream at $J$ and a left stream at $K+1$. 
Note that on the one hand,
$$  
\frac{C}{4} \left(\frac{1}{\gamma}\right)^{K-J}
\ge \frac {C}{4} > 
 2 \vareps D.
$$
Note also on the other hand that 
$$ \frac {C}{4} \left( \frac {1}{\gamma} \right)^{K-J} 
\le \frac {D}{4} \gamma^{L(J+1) + J - K } \le \frac {D}{4} \gamma^{LJ} \le \frac {\ol n - \ul n}{2}.$$
From Lemma \ref{lemma:conf} it therefore follows that $Y$ is confined to $[J,K+1]$ for at least
$$\frac{C}{4}  (1/ \gamma)^{K-J}$$ steps and that if, we let $\chn$ smallest integer such that
$$
\chn \ge \ul n + \frac{C}{4} \left(\frac{1}{\gamma}\right)^{K-J},
$$
then 
$$\chn \le \ol n.$$
Then, define
\begin{align*}
l_i
&:=
\ell(\chn, J+i)-\ell(\ul n, J+i),
&&
i\in\{1,\dots,K-J+1\},
\\
\label{didef1}
d_i
&:=
\Delta(\chn, J+i)-\Delta(\ul n, J+i),
&&
i\in\{1,\dots,K-J\}.
\end{align*}
Then, for $ l = (l_1, \ldots, l_{K-J +1})$, we see that
\begin{align*}
\norm{l}
=
\chn - \ul n
\ge
\frac{C}{4} \left(\frac{1}{\gamma}\right)^{K-J}.
\end{align*}
By the inductive assumption, we have
$$
\max_{n\in [\ul n, \chn]} \max_{j\in\{J+1,\dots,K\}} \abs{\Delta(n,j)}
\le
C \left(\frac{1}{\gamma}\right)^{K-(J+1)}
\le
4 \gamma \norm{l}.
$$
and thus, $$\abs{d_j} \le  8 \gamma \norm{l}$$  for all $
j\in\{1,\dots,K-J\}$.
Our choice of $\gamma$ therefore ensures that
$$
d_0
=
\Delta(\chn, J)-\Delta(\ul n, J)<0
$$
which is in contradiction with the definition of $\ul n$ and the fact that $\chn \in [\ul n, \ol n ]$. Hence, we conclude that \eqref{bound_beginning_2} holds for all $Y \in \cB_+$.

\medbreak
We now want to derive the lower bound on the entire time-interval $[\ul n , \ol n]$.
We now let
\begin {align*}
&
\wh n:=\inf\left\{n\ge\ul n: \Delta(n,K+1)<-C \left(\frac{1}{\gamma}\right)^{K-J}\right\}
\\
&
\wt n:=  \max\left\{n\le\wh n: \Delta(n,K+1) \ge -\frac{C}{2} \left(\frac{1}{\gamma}\right)^{K-J} \right\}.
\end{align*}
Assume that
$
\wh n \le \ol n$. 
From \eqref{bound_beginning_2}, we know that
$
\wt n>\ul n$ for all $Y \in \cB_+$.
But
\begin{align*}
\max_{n\in[\wt n, \wh n]} \Delta(n,K+1)  & = \Delta ( \wt n , K+1 ) \le  -\frac{D}{2} \gamma^{(J+1)L} \left(\frac{1}{\gamma}\right)^{K-J}   + 1 
\\
&<
- \frac {D}{2} \gamma^{L \times L}  + 1 \le -3 \vareps D + 1 \le   - 2\vareps D,
\end{align*}
so that it follows that on the whole time span $[\wt n, \wh n]$ the walk is confined to the interval $[J, K+1]$. We define now
\begin{align*}
l_i
&:=
\ell(\wh n, J+i)-\ell(\wt n, J+i),
&&
i\in\{1,\dots,K-J+1\},
\\
d_i
&:=
\Delta(\wh n, J+i)-\Delta(\wt n, J+i),
&&
i\in\{1,\dots,K-J\}.
\end{align*}
Then, if $l = ( l_1, \ldots, l_{K-J+1})$, we get that
$$
\norm{l}
=
\wh n - \wt n
\ge
\frac{C}{4} \left(\frac{1}{\gamma}\right)^{K-J} .
$$
By the inductive assumption we have
$$
\max_{n\in [\ul n, \wh n]} \max_{j\in\{J+1,\dots,K\}} \abs{\Delta(n,j)}
\le
C \left(\frac{1}{\gamma}\right)^{K-(J+1)}
\le
4 \gamma  \norm{l}.
$$
and thus,
\begin{align*}
\abs{d_j} \le \gamma \norm{l},
\qquad
i\in\{1,\dots,K-J\}.
\end{align*}
Hence, our choice for $\gamma$ ensures that
\begin{align*}
d_{K-J+1}
=
\Delta(\wh n, K+1)-\Delta(\wt n, K+1)>0
\end{align*}
which is in conflict with the definition of $\wt n$ and $\wh n$. We conclude that \eqref{remains} holds for $j=K+1$, which concludes the proof of the Lemma.

\section {Concluding remarks}

We now make some remarks on the proof, and list a few open problems that are directly related to the models that we have investigated in the present paper. Other related problems are discussed in \cite {erschler_toth_werner_10}.

\begin {itemize}
 \item
Just as in many other self-interacting random walks, it is rather hard to get direct information on the dynamics of the walker. The strategy of the proof presented in the present paper is to use some a priori information about the local time profile, in order to deduce the properties of the walker. As a consequence, there are many intuitive results that can not derived in this way (and it would be of course very nice to prove them).
\item
A  natural guess is that in the case where $\alpha < \alpha_{L}$ and the confined  walk $Y$ actually visits all sites of the interval infinitely often, it should be the case that, if one defines
$$ 
\Lambda_n = ( \ell (n, 1), \ldots, \ell (n, L+1)) - ( n u_1, \ldots, n u_{L+1}),
$$
then the Markov chain $( Y_n ,\Lambda_n)_{n \ge 0}$ should be positive recurrent and have an invariant distribution. In particular, this would imply not only the convergence of $(\ell (n, 1), \ldots, \ell (n,L+1)) / n $ towards $(u_1, \ldots, u_{L+1})$ but it would give a much finer description of the limiting behavior of the profile.

Unfortunately, it seems quite difficult to get an explicit expression for such a stationary distribution. It would be actually sufficient to prove the tightness of some well-chosen functional in order to prove the existence of the stationary measure, but we have not been able to find a simple way to approach this.

\item
One very natural question that the previous approach would enable to tackle (but that could possibly be studied by other means too) is to show that for the cases where we have proved that  $\# {\mathcal L} = L+2 $ with positive probability, then the probability that ${\mathcal L}= K+2$ for some $K > L$ is zero.

\item
A related question deals of course with the behavior of the walk when $\alpha$ is equal to one of the critical values. It seems intuitively clear that when $\alpha= \alpha_{L+1}$, then the walk will almost surely not be stuck on $L+2$ sites. This is because of the fact that the corresponding $d_0$ and $d_{L+1}$ are equal to zero, so that for infinitely many times $n$, if $X$ would stick to $Y$ forever, it would be at the edge of the interval and have a probability bounded from below to jump out of it. But our proof (that controls only the first-order behavior of the local-time profile) is not able to control this, nor to prove that $\# {\mathcal L}$ is equal to $L+3$ with positive probability.

\item
Even though it is a trivial observation, we would like to notice that for generic $\alpha$'s, it is the case that there exist many $K > L+1$ such that one can find many non-negative
$(l_1, \ldots, l_{K+1})$ with $d_1 = \cdots = d_{K}=0$, $d_0 > 0$ and $d_{K+1}< 0$.
This happens in fact as soon as the distance between $(K+1) \theta$ and $2 \pi \Z$ is smaller than all distances between the $j \theta$'s and $2 \pi \Z$ for all $j  \in \{ 1, \ldots, K \}$ (so that for non-rational $\pi / \theta$, there are infinitely many such $K$'s).

Note that if one wishes to find a stationary measure for the couple $(Y, \Lambda) $ as described above, then one would need to exclude these larger values of $K$, which is an indication that it could in fact be a rather complicated issue.

\item
The arguments that we have developed in the present paper seem to be adaptable in order to derive analogous results in some cases where the self-interaction depends on the local time at more than the four edges neighbouring $X_n-1, X_n, X_{n+1}$ and $X_{n+2}$ from the left.
For instance, one could replace the definition of $\Delta (n, j )$
by
$$
 \alpha_k \ell (n, j_k) + \ldots + \alpha_1 \ell(n,j-1)+\alpha_0 \ell(n,j)- \alpha_0 \ell(n,j+1) - \ldots - \alpha_k \ell (j+1+k)
,$$
the needed condition for the arguments to go through is dealing with the roots of the corresponding polynomial i.e. with the behavior of the corresponding generalized Fibonacci sequence.
\end {itemize}

\medbreak

\noindent{\bf Acknowledgements.}
BT thanks the kind hospitality of Ecole Normale Sup\'erieure, Paris, where part of this work was done. The research of BT is partially supported by the Hungarian National Research Fund, grant no. K60708. The research of WW is supported in part by Research supported in part by ANR-06-BLAN-00058. The cooperation of the authors is facilitated by the French-Hungarian bilateral mobility grant Balaton/02/2008.

\bigskip
\noindent
Affiliations and e-mails of authors:
\\[10pt]
{\sc Anna Erschler},
CNRS, D\'epartement de Math\'ematiques, Universit\'e Paris Sud Orsay,\\
email: {\tt Anna.Erschler@math.u-psud.fr}
\\[10pt]
{\sc B\'alint T\'oth},
Institute of Mathematics, Budapest University of Technology,\\
email: {\tt balint@math.bme.hu}
\\[10pt]
{\sc Wendelin Werner},
D\'epartement de Math\'ematiques, Universit\'e Paris Sud Orsay, and DMA, Ecole Normale Sup\'erieure\\
email: {\tt Wendelin.Werner@math.u-psud.fr}


\begin{thebibliography}{99}

\bibitem{amit_parisi_peliti_83}
D.\ Amit, G.\ Parisi, L.\ Peliti:
Asymptotic behavior of the `true' self-avoiding walk.
{\sl Phys.\ Rev.\ B}, {\bf 27}: 1635--1645 (1983)


\bibitem{erschler_toth_werner_10}
A. Erschler, B. T\'oth, W. Werner:
Some locally self-interacting walks on the integers,
preprint (2010).

\bibitem{obukhov_peliti_83}
S.\ P.\ Obukhov, L.\ Peliti:
Renormalisation of the ``true'' self-avoiding walk.
{\sl J.\ Phys.\ A}, {\bf 16}: L147--L151 (1983)


\bibitem{peliti_pietronero_87}
L.\ Peliti, L.\ Pietronero:
Random walks with memory.
{\sl Riv.\ Nuovo Cimento}, {\bf 10}: 1--33 (1987)


\bibitem{pemantle_92}
R. Pemantle:
Vertex-reinforced random walk.
{\sl Probab. Theory Rel. Fields}, {\bf 92}: 117-136 (1992)

\bibitem{pemantle_volkov_99}
R. Pemantle, S. Volkov:
Vertex-reinforced random walk on $\Z$ has finite range.
{\sl Ann. Probab.}, {\bf 27}: 1368-1388 (1999)

\bibitem{tarres_04}
P. Tarr\`es:
Vertex-reinforced random walk on $\Z$ eventually gets stuck on five points.
{\sl Ann. Probab.}, {\bf 32}: 2650-2701 (2004)


\bibitem{ttv}
{
P. Tarr\`es, B. T\'oth, B. Valk\'o (2010),
Diffusivity bounds for 1d Brownian polymers.
to appear in Ann. Probab., http://arxiv.org/abs/0911.2356.
}

\bibitem{toth_95}
B. T\'oth:
`True' self-avoiding walk with bond repulsion on $\Z$: limit theorems.
{\sl Ann. Probab.}, {\bf 23}: 1523-1556 (1995)

\bibitem{toth_01}
B. T\'oth:
Self-interacting random motions.
In: {\sl Proceedings of the 3rd European Congress of Mathematics}, Barcelona 2000, vol. 1, pp. 555-565, Birkhauser, 2001.

\bibitem{toth_werner_98}
B. T\'oth, W. Werner:
The true self-repelling motion.
{\sl Probab. Theory Rel. Fields}, {\bf 111}: 375-452 (1998)

\end{thebibliography}
\end{document}